\documentclass[12pt]{article}
\usepackage{group21,epsfig}

\begin{document}
\ifx\href\undefined\else\errmessage{Don't use HyperTeX!}\fi

\title{Path integrals and low-dimensional topology}
\author{Bogus\l aw Broda}
\address{Department of Theoretical Physics, University of
\L\'od\'z,\\ Pomorska 149/153, PL--90-236 \L\'od\'z, Poland\\
e-mail: \tt bobroda@krysia.uni.lodz.pl}

\maketitle

\section{Introduction}
The aim of our talk is to present a specific, non-perturbative, path-integral
approach
to topological invariants of knots/links and manifolds of dimension three and
four.
The technique is not rigorous but very intuitive and strongly motivated by
physics. An
exception is the four-dimensional case in Sect.~4, which is rather rigorous but
less
intuitive from physical point of view.

The plan of the paper is as follows. In Sect.~2, we will give an account of a
non-perturbative, path-integral derivation of standard ({\it quantum-group})
knot and
link invariants in the spirit of an original idea of Witten. We will also
present two
higher-dimensional generalisations. Sect.~3 is devoted to a Chern-Simons
approach to
topological invariants of three-dimensional manifolds of Reshetikhin, Turaev
and
Witten. Sect.~4 presents a four-dimensional version of the three-manifold
invariant.

\section{Knot and link invariants}
According to Witten \cite{W} topological invariants of knots and links can be
described as expectation values of Wilson loop observables in the framework of
three-dimensional Chern-Simons gauge theory. 

\subsection{Chern-Simons theory}
The classical action of Chern-Simons theory is given by an appropriately
normalized
{\it secondary characteristic class}
$$
S_{\rm CS}(A)={1\over4\pi}\int{\rm Tr}\left(A\wedge dA+{2\over3}A\wedge A\wedge
A\right),
\eqno{(1)}
$$
where $A$ is the gauge-potential (connection) for the simple, compact Lie group
$G$.
Since a knot (link) is a simple loop (collection of loops), we can associate an
observable to it: {\it Wilson loop(s)}. The Wilson loop is defined as
$$
W^l_{\cal C}(A)={\rm Tr}_l {\rm P}\exp\oint_{\cal C}A,
\eqno{(2)}
$$
where $\cal C$ is a loop (knot) in three-dimensional space, and $l$ labels
irrep's of
$G$. Since the action and Wilson loops are gauge-invariant and metric
independent, we
can claim that expectation values
$\left<W_{\cal C}(A)\right>=
\int W_{\cal C}(A)\exp\left[ikS_{\rm CS}(A)\right] DA$
should be metric independent as well, i.~e. they are topological invariants.
$k$ is an
integer (coupling constant or conformal weight). To determine the corresponding
topological invariant, one should derive the, so-called, {\it skein relation},
allowing to recursively compute the topological invariant. An example of a
three-component skein relation is given below
$$
\alpha
\setbox0=\hbox{
\setbox1=\hbox{\bigg/}
\setbox2=\hbox{$\backslash$}
\hbox to -\wd1{} \copy1
\kern-\wd1\raise.5\ht1\copy2
\lower\dp1\hbox{\raise\dp2\copy2}
}\copy0
+\beta
\setbox0=\hbox{
\setbox1=\hbox{$\bigg\backslash$}
\setbox2=\hbox{/}
\hbox to -\wd1{} \copy1
\kern-\wd1\lower.5\ht1\copy2
\raise\dp1\hbox{\lower\dp2\copy2}
}\copy0
+\gamma
\;\bigg\vert\;\bigg\vert\;
=0.
$$
To perform a concrete calculation we translate the operator form of the Wilson
loop
(2) to a path-integral one, and apply the Stokes theorem transforming a line
integral
into a surface one \cite{B1}. If one of the lines pierces the surface, there is
a
contribution to the path integral we should calculate \cite{B2}. As a result,
we
obtain all well-known link invariants \cite{B3}:

\bigskip

\begin{tabular}[]{c|c|c|}
                  &$\quad$Fundamental
                   representation$\quad$
                                         &$\quad$Higher
                                         representations$\quad$       \\
\hrulefill        &\hrulefill            &\hrulefill                     \\
$\quad SU(2)\quad$
                  &Jones, Kauffman       & Akutsu-Wadati                 \\
   $SU(n)$        & HOMFLY               & \dots                         \\
   $SO(n)$        & Dubrovnik Kauffman   & \dots                         \\
\end{tabular}

\subsection{Higher dimensions}
It is not difficult, at least formally, to generalise our considerations to the
case
of higher dimensions. In principle, there are two possibilities we will very
shortly
describe.

\subsubsection{Inhomogeneous Chern-Simons theory}
One can consistently define, in arbitrary dimension, an action of the form (1),
with
$A=\sum_{i=\rm odd}A_i$ inhomogeneous form, where $A_i$ is a $i$-form.
Corresponding
generalisation of observables is straightforward \cite{B4}. This way one can
describe
invariants of links analogous to the three-dimensional case.

\subsubsection{BF-theory}
In $d$ dimensions the, so-called, {\it BF}-theory is defined by
$S_{\rm BF}(A,B)=\int{\rm Tr} B\wedge F,$
where $F$ is the gauge strength (curvature), and $B$ is an independent
non-abelian
($d-2$)-form field. This theory naturally describes linking phenomena between
components of dimension 1 and $d-2$ \cite{B5}.

\subsubsection{Quantisation}
There is a purely technical but important issue concerning the quantisation
procedure.
In higher dimensions, theories under consideration are plagued by on-shell
reducible
gauge symmetries. To cope with this problem one should use the formalism of
Batalin
and Vilkovisky. This procedure introduces a host of different kinds of ghosts
and
auxiliary fields (and also metric) \cite{B6}. 

\section{Three-manifold invariants}
Up to the present moment, we have been interested in knot and link invariants
of knots
and links in ${\cal S}^3$ (or $R^3$). It is time now, to extend our analysis to
the
case of an arbitrary closed, connected, three-dimensional manifold ${\cal
M}^3$. First
of all, we will be interested in the extreme case of empty knot/link, i.~e.\ in
the
manifold ${\cal M}^3$ itself. As far as {\it non-perturbative calculations} are
concerned, there is a class of basically equivalent invariants which
definitions
depend on topological description of ${\cal M}^3$. There are the three main
possibilities: (1) surgery on a link, (2) Heegaard decomposition, (3)
simplicial
decomposition (triangulation). Correspondingly, we have three types of
invariants: (1)
Reshetikhin-Turaev-Witten (RTW), (2) Kohno, (3) Turaev-Viro. In this talk, we
will
confine ourselves to the first type.

\subsection{RTW invariant}
The idea is to use the fact that it is possible to obtain an arbitrary
three-dimensional manifold ${\cal M}^3$ via surgery on a link, i.~e. by
attaching
two-handles along a link $L$. In physical terms, {\it attaching} means gluing
or
identifying boundary values \cite{B7}.

\subsubsection{Second Kirby move}
In the first step, we should formally calculate the partition function of
Chern-Simons
theory with fixed boundary conditions (holonomies) along all components of the
surgery
link $L$
$$
Z(g_1,\dots,g_N)=
\left< \delta(g_1,{\rm Hol}_{{\cal C}_1})
\cdots
\delta(g_N,{\rm Hol}_{{\cal C}_N})
\right>,
$$
where $\delta$ is a (group-theoretic) Dirac delta-function,
$\delta(g,h)=
\sum_l \overline{\chi_l(g)} \chi_l(h),$
$\chi_l$ are characters, and ${\rm Hol}_{{\cal C}_i}$ are holonomies around
${\cal
C}_i$. Then
$$
Z(g_1,\dots,g_N)=
\left< \sum_{l_1} \overline{\chi_{l_1}(g_1)} W_{{\cal C}_1}^{l_1}
\cdots
\sum_{l_N} \overline{\chi_{l_N}(g_N)} W_{{\cal C}_N}^{l_N}
\right>.
$$
In the second step, we should identify and sum up corresponding boundary values
($g_1$, ... , $g_N$). Thus, we obtain
$$
Z_{{\cal M}^3}=\int dg_1 \cdots dg_N Z_O(g_1^{-1}) \cdots Z_O(g_N^{-1})
Z(g_1, \dots, g_N)
=\left< \omega_{{\cal C}_1} \cdots \omega_{{\cal C}_N} \right>,
\eqno{(3)}
$$
where $\omega_{{\cal C}_i}=
\sum_l \left<W_O^l(A)\right> W_{{\cal C}_i}^l(A),$
and ``O'' refers to calculations performed for an unknot.
The quantity (3) is invariant with respect to the, so-called, second Kirby
move.

\subsubsection{First Kirby move}
To stabilise (3) to obtain a true topological invariant we should normalize it
with
the factor (in the denominator)
$\left<\omega_+\right>^{b_+} \left<\omega_-\right>^{b_-},$
where $b_+$ ($b_-$) means the number of positive (negative) eigenvalues of the
liking
matrix ${\rm L}\!{\rm k}$ of $L$, and $+$ ($-$) as a subscript pertains to an
unknot
with positive (negative) twist, correspondingly. The normalized (3) is
independent of
the both Kirby moves. 

\subsection{Knots and links}
Incorporation of knots and links is straightforward. We should modify (3) to
the
expression
$$
\left< W^{l_1}_{{\bar{\cal C}}_1} \cdots W^{l_M}_{{\bar{\cal C}}_M}
\omega_{{\cal C}_1} \cdots \omega_{{\cal C}_N}
\right>,
$$
where {\it barred} quantities mean components of an ordinary knot/link.

\section{Generalisation to four dimensions}
We can easily generalise our considerations to the case of a four-dimensional
(closed,
connected) manifold ${\cal M}^4$. There is an analogous possibility to
construct an
arbitrary topological ${\cal M}^4$ via (generalised) surgery on a special link
$\ell$
\cite{Sa}. The special link $\ell$ consists of two kinds of components: {\it
undotted}
$L$ and {\it dotted} $\dot L$. The four-dimensional Kirby calculus is a little
bit
more complicated than in the three-dimensional case. There are three {\it
second}
Kirby moves now:

\begin{figure}[h]
\begin{center}
\mbox{
\epsfig{file=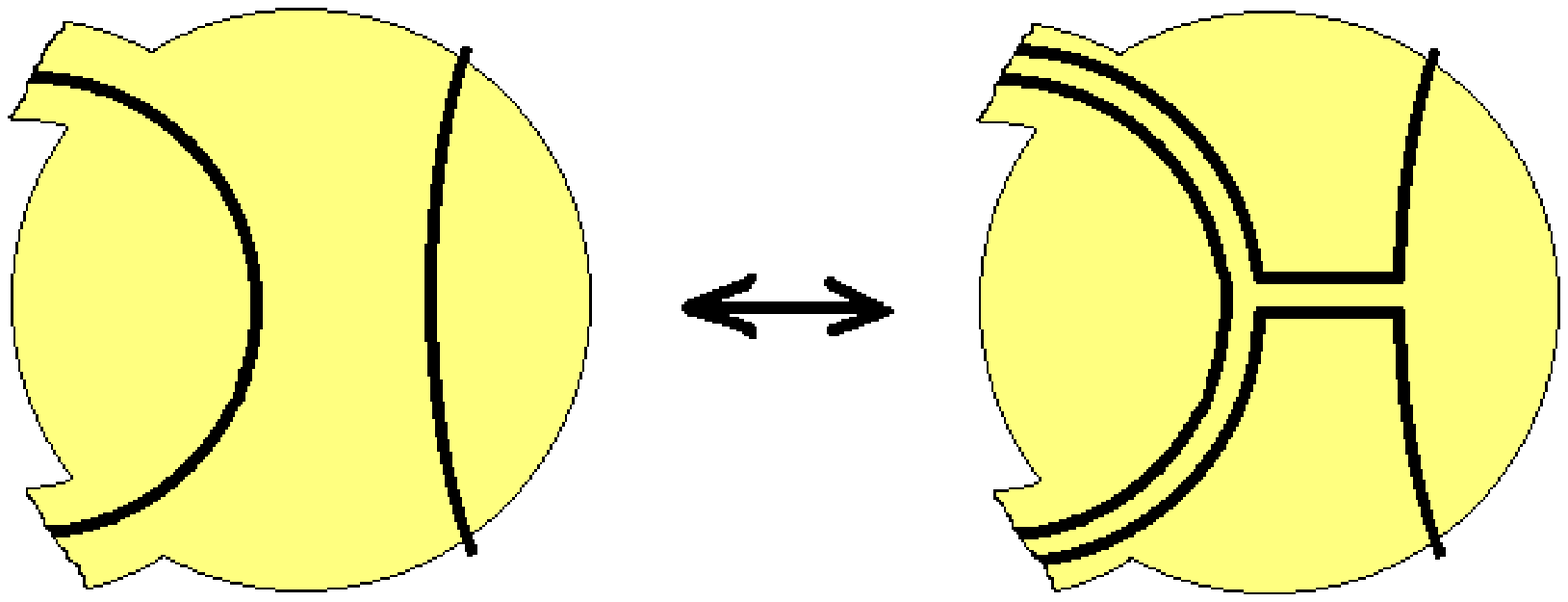, width=4.5cm}
\hfil
\epsfig{file=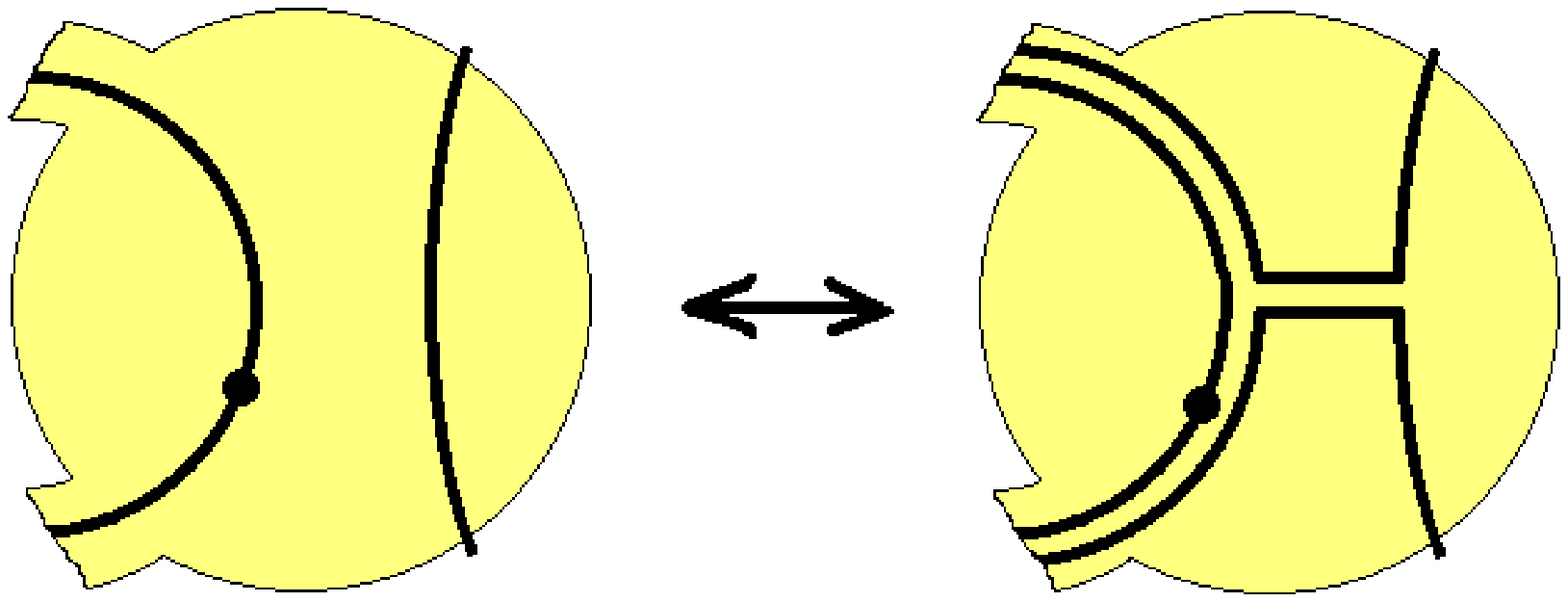, width=4.5cm}
\hfil
\epsfig{file=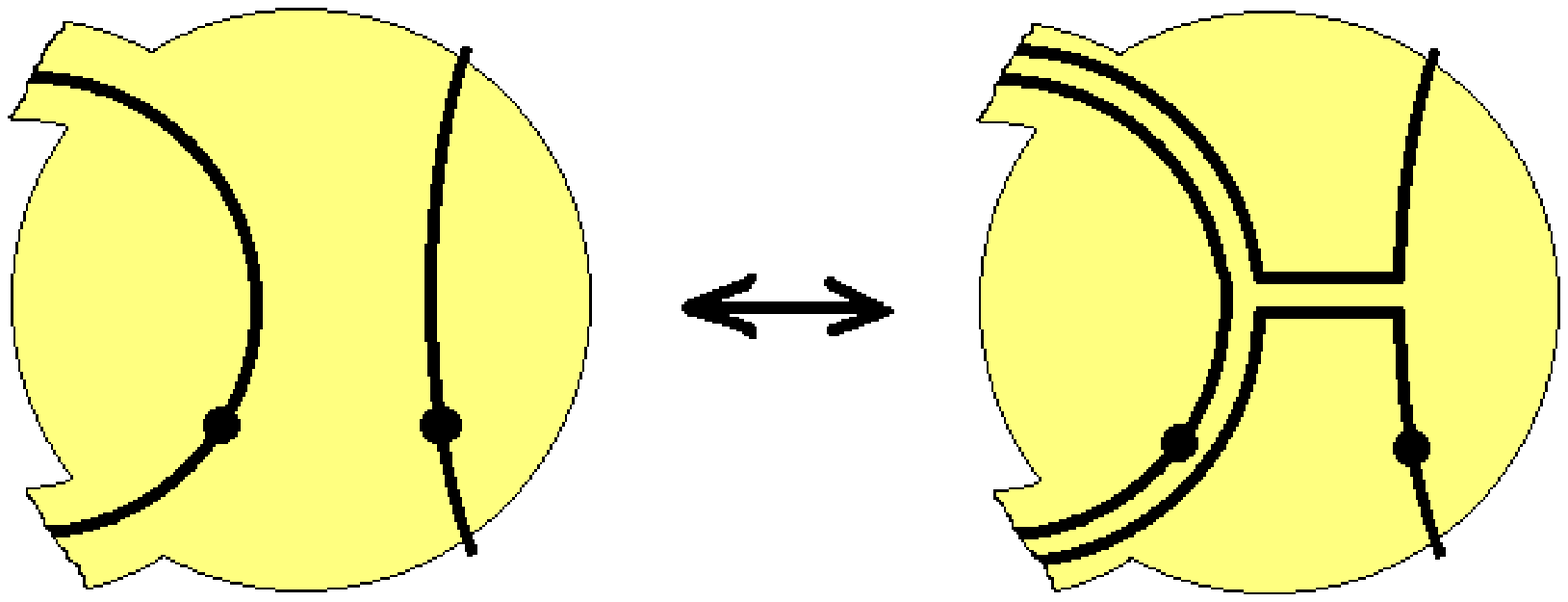, width=4.5cm}
}
\end{center}
\end{figure}

Nevertheless, we can find, purely combinatorically, a topological invariant
\cite{B8}
$$
{\cal I}_k ({\cal M}^4)=
{\left< \omega_{{\cal C}_1}^+ \cdots \omega_{{\cal C}_N}^+
\omega_{{\cal C}_{\dot 1}} \cdots \omega_{{\cal C}_{\dot N}} \right>
\over
\left< \omega_O^+ \right>^\nu
\left< \omega_H^+ \omega_{\dot H} \right>^{(N+\dot N - \nu)/2}},
\eqno{(4)}
$$
where for simplicity we have assumed that $G=SU(2)$. Here $\omega^+$ means {\it
even}
part of $\omega$ (only integer spins), $\nu$ is the nullity of the linking
matrix
${\rm L}\!{\rm k}$ of $\ell$, and $H$, $\dot H$ are components of the {\it Hopf
link}.

\section{Finishing remarks}
We are able to yield all invariants of knots, links and low-dimensional
manifolds (of
dimension three and four) pertaining to quantum groups (and their
representations)
using their classical counterparts. Other interesting chapter on topological
invariants could contain an approach to {\it perturbative} invariants like
Vassiliev
invariants of knots, Casson-Walker invariant, and its {\it higher order}
generalisations, and dimensionally reduced {\it Seiberg-Witten invariant}.

\section*{Acknowledgments}
I am greatly indebted to Prof.~H.-D.~Doebner for his kind hospitality during my
stay
in Clausthal.
I would also like to thank the organisers, and in particular
\hbox{Prof.~H.-D.~Doebner} for giving me the opportunity to present my results
during
the Colloquium.

The paper has been supported by Polish grant {\tt 2 P03B 094 10} and by the
Humboldt
Foundation.

\end{document}